\documentclass[a4paper,12pt]{article}

\usepackage{amsmath,amssymb,amsthm,latexsym,graphicx}
\usepackage{marvosym}
\usepackage{dictsym}
\usepackage{natbib}
\usepackage{longtable}
\usepackage{hyperref}

\usepackage{titling}
\usepackage{a4wide}
\usepackage[compact]{titlesec}
\usepackage{enumitem}
\setlist{nolistsep}

\hypersetup{
	colorlinks=true,
    breaklinks=true,
    filecolor=blue,
    urlcolor=blue,
    linkcolor=blue,
    citecolor=blue,
    anchorcolor=blue,
    frenchlinks=true,
    pdfborder={0 0 0},
    pdfdisplaydoctitle=true,
    bookmarksnumbered=true
}

\theoremstyle{definition}

\newtheorem{lemma}{Lemma}
\newtheorem{proposition}[lemma]{Proposition}
\newtheorem{theorem}[lemma]{Theorem}

\newtheorem{definition}[lemma]{Definition}
\newtheorem{example}[lemma]{Example}

\newtheorem{tab}[lemma]{Table}
\newtheorem{algorithm}[lemma]{Algorithm}

\title{Algorithms for singularities and real structures of weak Del Pezzo surfaces}
\author{Niels Lubbes}

\date{\today}

\begin{document}\maketitle
\begin{abstract}
 In this paper we consider the classification of singularities (Du Val) and real structures (Wall) of weak Del Pezzo surfaces from an algorithmic point of view.

 It is well known that the singularities of weak Del Pezzo surfaces correspond to root subsystems. We present an algorithm which computes the classification of these root subsystems. We represent equivalence classes of root subsystems by unique labels. These labels allow us to construct examples of weak Del Pezzo surfaces with the corresponding singularity configuration.

 Equivalence classes of real structures of weak Del Pezzo surfaces are also represented by root subsystems. We present an algorithm which computes the classification of real structures. This leads to an alternative proof of the known classification for Del Pezzo surfaces and extends this classification to singular weak Del Pezzo surfaces.

 As an application we classify families of real conics on cyclides.
\end{abstract}\begingroup \def\addvspace#1{}

\tableofcontents

 \endgroup

\section{Introduction}

\subsection{Problems}
Many properties of weak Del Pezzo surfaces are determined by its Picard group and the effective (-2)-classes therein. The (-2)-classes of a given weak Del Pezzo surface form a root system and the effective (-2)-classes form a root subsystem. These root subsystems are well studied, and the classification can for example be found in \cite{val1} and \cite{dol1}. In \cite{wal1} these root subsystems are used to study real Del Pezzo surfaces. Since this subject is studied quite extensively it is not easy to cite all the relevant literature, but a nice overview can be found in \cite{dol1}.

 In this paper we will address the following two classifications from an algorithmic point of view:
\begin{itemize}\addtolength{\itemsep}{1pt}
\item[$\bullet$] the classification of root subsystems of weak Del Pezzo surfaces (\cite{val1},\cite{dol1}), and
\item[$\bullet$] the real structures of weak Del Pezzo surfaces (\cite{wal1}).
\end{itemize} More precisely, in this paper we address to the following 2 problems.

 \textbf{Problem 1.} Compute a table which classifies weak Del Pezzo surfaces up to Weyl equivalence. Provide methods which, for each entry in this table, enables to do explicit calculations in the Picard group of weak Del Pezzo surfaces in the corresponding equivalence class. Moreover provide methods to compute explicit parametrizations of weak Del Pezzo surfaces in a given equivalence class.

 A solution to this problem is provided by Algorithm~\ref{alg:f1_geo_basis} and \textsection\ref{sec:f1_examples}.

 \textbf{Problem 2.} Compute a table which classifies conjugacy classes of real structures on weak Del Pezzo surfaces. Provide explicit coordinate representations for the action of real structures on the Picard group.

 A solution to this problem is provided by Theorem~\ref{thm:f1_real_alg}. The method of classifying the real structures does not use the existing classification in \cite{wal1} and thus provides an alternative proof. In \cite{wal1} only smooth Del Pezzo surfaces are considered. In \textsection\ref{sec:f1_cyclide} we consider real structures on weak Del Pezzo surfaces with a singular anticanonical model.

 As an application of the algorithms in this paper we classify families of real conics on cyclides. This gives an coordinate independent approach to the classification of cyclides in \cite{tak2}. Implementations of the algorithms in this paper can be found on the homepage of the author.

 \subsection{Motivation}
A \textit{minimal family}\index{minimal family} is defined as a 1-dimensional family of rational curves of minimal degree, that can be defined by the fibers of a morphism. In \cite{nls1} we reduce the classification of minimal families to the classification of minimal families, which generate either geometrically ruled surfaces or weak Del Pezzo surfaces. The unique minimal families of a geometrically ruled surface is defined by its ruling. For a weak Del Pezzo surface the minimal families (that can be defined by fibers of a morphism) are uniquely defined by special classes in its Picard group together with intersection product and Betti numbers. For this reason we say that weak Del Pezzo surfaces carry equivalent minimal fibration families  if and only if  they are Weyl equivalent as in problem 1. See \cite{nls-f3} for more information.

 From \cite{Schicho:00a} we know that surfaces that carry at least 2 families of conics are either weak Del Pezzo surfaces or geometrically ruled. This paper extends this result by also incorporating the real structure.

 Blum conjecture states that maximum number of ways a surface, which is not the sphere or the plane, can be generated by a family of circles is $6$ (see \cite{blum1}). This conjecture was proven by \cite{tak1} depending on differential geometric results. As an application of this paper we prove an alternative conjecture  with respect to  the Euclidean topology: a compact surface contains either infinite or at most 6 families of real conics. From \cite{Schicho:00a} it follows that we only have to look at Table~\ref{tab:f1_real} with rank $<7$. It can be shown that the f1 column denotes the number of real lines and thus must be 0 since we assume that the surface is compact. It follows by inspection that the maximum is attained at index 13. Strictly speaking the latter conjecture is not equivalent to Blum conjecture, but it is possible to prove Blum conjecture using the results of this paper. This will be addressed to in \cite{nls-f4}.

 Minimal families have applications in geometric modeling (see for example \cite{pot2}). Using our classification we can construct surfaces with required properties.

 The algorithms for problem 1 in this paper are improvements of the algorithms in the PhD thesis: \cite{nls2}. Since there are more than 200 cases to consider, our experience is that mistakes are easily made when computed by hand. I hope that this paper is of service to others and saves a lot of time.

\section{Weak Del Pezzo surfaces and root subsystems}

\label{sec:f1_def}Let ${\textrm{X}}$ be a nonsingular complex projective surface. The \textit{enhanced Picard group}\index{enhanced Picard group} $A({\textrm{X}})$ of ${\textrm{X}}$ is defined as $ (~\textrm{Pic}({\textrm{X}}),~K,~\cdot,~ h~) $ where $\textrm{Pic}({\textrm{X}})$ is the Picard group, $K$ is the canonical divisor class of ${\textrm{X}}$, $\ensuremath{\textrm{Pic}({\textrm{X}})\times\textrm{Pic}({\textrm{X}})\stackrel{\cdot}{\rightarrow}{\textbf{Z}}}$ is the intersection product on divisor classes, and $\ensuremath{{\textbf{Z}}\times\textrm{Pic}({\textrm{X}})\stackrel{h}{\rightarrow}{\textbf{Z}}_{\geq0}}$ assigns the $i$-th Betti number to a divisor class for $i\in{\textbf{Z}}$. For $h(i,D)$ we use the notation $h^i(D)$ and we denote $D\cdot C$ by $DC$  for $D,C\in \textrm{Pic}({\textrm{X}})$.

 We consider enhanced Picard groups isomorphic  if and only if  there exists an isomorpism of the Picard groups that preserves the canonical divisor class and is compatible with $\cdot$ and $h$.

 We call ${\textrm{X}}$ a \textit{weak Del Pezzo surface}\index{weak Del Pezzo surface}  if and only if  its anticanonical class $-K$ is nef and big. For a weak Del Pezzo surface with $K^2<8$ we have that $A({\textrm{X}})={\textbf{Z}}\langle ~H,~ Q_1,~ \ldots,~ Q_r~\rangle $ with $H^2=1$, $Q_iQ_j=-\delta_{ij}$, $HQ_i=0$, $K=-3H+Q_1+\ldots+Q_r$ and $r=9-K^2$. See \cite{dol1} for weak Del Pezzo surfaces.

 The \textit{(a,b)-set}\index{(a,b)-set} of $A({\textrm{X}})$ is defined as $\{-CK=a {\textrm{~and~}} CC=b|C\in A({\textrm{X}})\}$. Let $F({\textrm{X}})$ be the (0,-2)-set. The \textit{Weyl object}\index{Weyl object} $W({\textrm{X}})$ is defined as a tuple $(R, S)$ where $R=F({\textrm{X}})$ and $S=\{ \pm C |C\in F({\textrm{X}}) {\textrm{~and~}} h^0(C)\geq 0\}$.

 We call a root system a \textit{C1 root system}\index{C1 root system}  if and only if  its Dynkin type is either $A_1$, $A_1+A_2$, $A_4$, $D_5$, $E_6$, $E_7$ or $E_8$. If ${\textrm{X}}$ is a weak Del Pezzo surface then $R$ is a C1 root system in the vectorspace ${\textbf{R}}\langle C\in  \textrm{Pic}{\textrm{X}} ~~|~~ CK=0\rangle $ (proposition 8.2.10 in \cite{dol1}). Moreover, we have that $S\subset R$ is a root subsystem. See \cite{bou1} and \cite{gra1} for the theory on root systems.

 We recall that 2 root subsystems $S,S'\subset R$ are isomorphic  if and only if  there exists an action $\ensuremath{R\stackrel{w}{\rightarrow}R}$ of the Weyl group on $R$  such that  $w(S)=S'$.
\begin{proposition}\label{prop:f1_weyl}\textbf{\textrm{(properties of Weyl objects of weak Del Pezzo surfaces)}}
 Let ${\textrm{X}}$ and ${\textrm{X}}'$ be weak Del Pezzo surfaces. Let $A({\textrm{X}})$ be the enhanced Picard group. Let $W({\textrm{X}})$ be the Weyl object (thus a root subsystem).

We have that ${\textit{A}}({\textrm{X}})\cong{\textit{A}}({\textrm{X}}')$  if and only if  ${\textit{W}}({\textrm{X}})\cong{\textit{W}}({\textrm{X}}')$.

\begin{proof}
 See section 8.2.8 and corollary 8.2.33 in \cite{dol1} or appendix E, section 4 in \cite{nls2}.  \end{proof}
\end{proposition}
In order to algorithmically check whether root subsystems are isomorphic we introduce the following definition and proposition.
\begin{definition}
\label{def:}
\textrm{\textbf{(complement and double complement of root subsystems)}}
  Let $S\subset R$ be a root subsystem. The \textit{complement root subsystem}\index{root system!complement root subsystem}\index{complement root subsystem} is defined as $C_R(S)=\{ r\in R | (r,s)=0 \textrm{ for all } s\in S \}$. The \textit{double complement root subsystem}\index{root system!double complement root subsystem}\index{double complement root subsystem} is defined as $K_R(S)=\{ r\in R | r \in {\textbf{R}}\langle S\rangle  \}$.

\end{definition}

\begin{proposition}\label{prop:f1_C_K}\textbf{\textrm{(properties of complement and double complement of root subsystems)}}

 Let $S$ and $S'$ be root subsystems of an irreducible C1 root system $R$.
 Let $D(S)$ denote the Dynkin diagram of $S$.

 We have that $(S,R)\cong(S',R)$  if and only if  $D(S)\cong D(S')$ (as graphs) , $|C_R(S)| =  |C_R(S')|$ and $|K_R(S)| = |K_R(S')|$.

\begin{proof}

\textit{Claim:} We have that $\Longrightarrow$.
\\  We have that $C_R(S)={\textbf{R}}\langle S\rangle ^{\bot}\cap R$ where ${\textbf{R}}\langle S\rangle  \subset V$ is an Euclidean vector subspace.  Similarly we have that $K_R(S)={\textbf{R}}\langle S\rangle ^{\bot\bot}\cap R$.  We have an action $\ensuremath{R\stackrel{f}{\rightarrow}R}$ of an element the Weyl group  such that  $f(S)=S'$.  This extends linearly to an Euclidean isometry $\ensuremath{V\stackrel{f'}{\rightarrow}V}$  such that  $f'(R)=R$.  From $f'$ preserving orthogonality it follows that this claim holds.

\textit{Claim:} We have that $\Longleftarrow$.
\\  This claim follows from tables 10.1 and 10.2 in \cite{osh1}.  These tables contain $D(C_R(S))$ for each root subsystem $S$ of C1 root systems.  We find that if $D(S)\cong D(S')$, $C_R(S)\cong  C_R(S')$ then $(S,R)\cong(S',R)$, for almost all $S,S'$.  Except for the subsystems $4A_1$, $A_5+A_1$ and $A_3+2A_1$ in $E_7$ we see that there are two Weyl equivalence classes which cannot be distinguished this way.  For example there exists root subsystems $S$ and $S'$ of root system $R$  such that  $D(R)=E_7$,  $D(S)=D(S')=A_5+A_1$ and $D(C_R(S))=D(C_R(S'))=\emptyset$.  However for these cases we find that $K_R(S)\ncong_b K_R(S')$ and thus this claim holds.
 \end{proof}
\end{proposition}

\section{Algorithm for C1 root subsystems classification}
Motivated by Proposition~\ref{prop:f1_weyl} we would like to represent isomorphism classes of Weyl objects of weak Del Pezzo surfaces by a label. We would like to have `geometric' labels which allow us to construct examples and to do explicit computations.
\begin{definition}
\label{def:f1_C1_labels}
\textrm{\textbf{(C1 label)}}
 The set of \textit{C1 label elements}\index{C1 label elements} is defined as $G=\{$ $\pm ab$, $\pm 1abc$,$\pm 2ab$, $\pm 30a$ $~~|~~$ $a,b,c\in[1, \ldots, 8]$ $\}$, A \textit{C1 label}\index{C1 label} is defined as $(L,r)$ where $L \subset G$ and $r\in [1, \ldots, 8]$. Let ${\textbf{Z}}\langle ~H,~ Q_1,~ \ldots,~ Q_r~\rangle $ be a basis. The \textit{C1 label elements to root function}\index{C1 label elements to root function} $\ensuremath{G\stackrel{\varphi}{\rightarrow}R}$ is defined as follows:
\begin{itemize}\addtolength{\itemsep}{1pt}
\item[$\bullet$] $\varphi(ab)    =   Q_a-Q_b$,
\item[$\bullet$] $\varphi(1abc)  =   H-Q_a-Q_b-Q_c$,
\item[$\bullet$] $\varphi(2ab)   =   2H-\ensuremath{\overset{}{\underset{i\in[1,\ldots, 8],i\notin\{a,b\}}{\sum}}}Q_i$, and
\item[$\bullet$] $\varphi(30a)   =   3H-\ensuremath{\overset{}{\underset{i\in[1,\ldots, 8],i\notin\{a\}}{\sum}}}Q_i-2Q_a$,
\end{itemize} for all $a,b,c\in[1,\ldots, 8]$ and $\varphi(-l)=-\varphi(l)$ for all $l\in G$. We call a C1 label \textit{geometric}\index{C1 label!geometric}\index{geometric}  if and only if  the label does not contain minus symbols.

\end{definition}
We will use C1 labels $(L,r)$ for representing a basis $(~\varphi(l)~)_{l\in L}$ of a root subsystem $S$ of a C1 root system $R$. In other words we can represent isomorphism classes of root subsystems by C1 labels.

 The following algorithm computes geometric C1 labels for C1 root subsystems. Later we will use the output to construct explicit examples, and to do computations in the Picard group of weak Del Pezzo surfaces. The classification of C1 root subsystems is well known (see section 8.2.3 in \cite{dol1}).
\begin{algorithm}\label{alg:f1_geo_basis}\textbf{\textrm{(geometric bases of C1 root subsystems)}}{\small

 \textbf{Output:} A list of geometric C1 labels. Each geometric C1 label represents a basis for a unique representative of each isomorphism class of root subsystems.
\begin{itemize}\addtolength{\itemsep}{1pt}
\item[$\bullet$] Let $\ensuremath{G\stackrel{\varphi}{\rightarrow}R}$ be the C1 label elements to root function. By abuse of notation we will not distinguish between a C1 label $(L,r)$, and the corresponding root set $(~\varphi(l)~)_{l\in L}$.
\item[$\bullet$] Let \\
$Z_0 := $ {\tiny$\{$ $12$, $23$, $34$, $45$, $56$, $67$, $78$, $1123$, $-1145$, $-1345$, $-1167$, $-1178$, $-278$, $-218$, $-308$, $218$ $\}$}, \\
$Z_1 := $ {\tiny$\{$ $1123$, $12$, $23$, $34$, $45$, $56$, $67$, $78$ $\}$}, \\
$Z_2 := $ {\tiny$\{$ $1123$, $12$, $23$, $34$, $45$, $56$, $67$, $78$, $1145$, $1347$, $1678$, $1127$, $1456$, $1567$, $234$, $278$, $308$ $\}$}, \\
$Z_3 := $ {\tiny$\{$ $1123$, $1345$, $1165$, $1285$, $1673$, $1274$, $1684$, $1178$ $\}$}.
\item[$\bullet$] We apply the following two steps by looping over $r \in [2,\ldots,8]$ and $i\in[0,\ldots,3]$.
\item[$\bullet$] Let $Y_i$ be the set of C1 labels $(L,r)$ where $L$ is a subset of $Z_i\cap\{\pm ab, \pm 1abc,\pm 2ab, \pm 30a|a,b,c\in[r]\}$  such that  $L$ is linear independent and $(l,m)\in \{0,1\}$ for all $l,m \in L$. So $L$ defines a valid basis for a root subsystem of rank $r$.
\item[$\bullet$] Using Proposition~\ref{prop:f1_C_K} we reduce the set $Y_i$  such that  it only contains non-isomorphic root subsystems represented by C1 labels $(L,r)$. Sage contains a function which checks isomorphisms of graphs. We check for the following special case: the C1 labels $(\{1123\},3)$ and $(\{23\},3)$ are non-isomorphic root systems of Dynkin type $A_1$ in the reducible (!) root system of type $A_1+A_2$ (although their complement and double complement root subsystems are isomorphic).
\item[$\bullet$] We have that $\cup_i Y_i$ contains isomorphic root subsystems. For a given equivalence class of root subsystems in $Y_0$ the preferred representative is in $Y_1$. If it does not exists in $Y_1$ then $Y_2$. If it does not exists in $Y_2$ then $Y_3$. If it does not exists in $Y_3$ then $Y_0$.
\item[$\bullet$] We return the obtained set of C1 labels.
\end{itemize} }
\end{algorithm}
The output of Algorithm~\ref{alg:f1_geo_basis} is in Table~\ref{tab:f1_C1}.
\begin{proposition}\label{prop:}\textbf{\textrm{(geometric bases of C1 root subsystems)}}

  Algorithm~\ref{alg:f1_geo_basis} is correct.

\begin{proof}
 The set $Z_0$ is constructed using Dynkin's Algorithm, which recursively adds the negative of maximal roots of sub diagrams. This ensures independently of existing classifications that all isomorphism classes of root subsystems is reached in $Z_0$.
 The set $Z_i$ for $i>0$ are constructed  such that  each equivalence class of root subsystems has a geometric C1 label. The set $Z_1$ forms a basis for the root system of Dynkin type $E_8$. If we apply the algorithm using only $Z_0$ and $Z_1$ then the output contains C1 labels that are not geometric. We construct sets $Z_2$ and $Z_3$  such that  each of these non-geometric C1 labels is represented by a geometric C1 label. The root subsystem $4A_2$ in $E_8$ requires extra attention (see Example~\ref{ex:f1_exm_degree_one}) and leads to the labels of $Z_3$.
 From Proposition~\ref{prop:f1_C_K} it follows that no two C1 labels in Table~\ref{tab:f1_C1} represent isomorphic root subsystems of an irreducible root system.
 When the rank is 3 then we consider root subsystems of a reducible root system $A_1+A_2$.
 We have to make an explicit check for two non-isomorphic root subsystems of Dynkin type $A_1$.  \end{proof}
\end{proposition}

\section{Constructing examples of weak Del Pezzo surfaces}

\label{sec:f1_examples}For a given geometric C1 label $(L,r)$ in Table~\ref{tab:f1_C1} we want to construct examples of a weak Del Pezzo surface ${\textrm{X}}$ in the corresponding Weyl equivalence class. See chapter 7, section 4, remark 98 and remark 99 in \cite{nls2} for an alternative method for constructing weak Del Pezzo surfaces of degree one and two.

 Let $\ensuremath{G\stackrel{\varphi}{\rightarrow}R}$ be the C1 label elements to root function. We choose a set of points $(p_i)_{i\in[r]}$ in the projective plane that are generic except:
\begin{itemize}\addtolength{\itemsep}{1pt}
\item[$\bullet$] $p_{i_2}$ is infinitely near to $p_{i_1}$  if and only if  $\l\in L$  such that  $\varphi(l)=Q_{i_1}-Q_{i_2}$,
\item[$\bullet$] $p_{j_t}$ lie on a line  if and only if  $\l\in L$  such that  $\varphi(l)=H - Q_{j_1} - Q_{j_2} - Q_{j_3}$,
\item[$\bullet$] $p_{k_t}$ lie on a conic  if and only if  $\l\in L$  such that  $\varphi(l)=2H -  Q_{k_1} - Q_{k_2} - Q_{k_3} - Q_{k_4} - Q_{k_5} - Q_{k_6}$, or
\item[$\bullet$] $p_{l_t}$ lie on a cubic with a double point at $p_{l_1}$  if and only if  $\l\in L$  such that  $\varphi(l)=3H - 2Q_{l_1} - Q_{l_2} - Q_{l_3} - Q_{l_4} - Q_{l_5} - Q_{l_6} - Q_{l_7} - Q_{l_8}$,
\end{itemize} for $i_t, j_t, k_t, l_t \in [8]$. Note that for this construction we require that the C1 label $(L,r)$ is geometric.

 We construct a linear series of degree 3 polynomials with base locus $(p_i)_{i\in[r]}$ (see chapter 8 in \cite{nls2}). Let $\ensuremath{{\textbf{P}}^2\stackrel{f}{\dashrightarrow}{\textrm{Y}}}$ be the associated map of the linear series. We define $\ensuremath{{\textrm{X}}\stackrel{g}{\rightarrow}{\textbf{P}}^2}$ to be the blowup of the projective plane at the base locus of the linear series. We define $D\in\textrm{Pic} {\textrm{X}}$ to be the pullback of the hyperplane sections along $g$. We have that ${\textrm{Y}}=\varphi_D({\textrm{X}})$ with $D=-K$ being the anticanonical divisor class of ${\textrm{X}}$.
\begin{center}
 \begin{tabular}
{l@{~~}l@{}l} ${\textrm{X}}$ & & \\
 $g\downarrow$ & $\searrow \varphi_{D}$ & \\
 ${\textbf{P}}^2$  & $\stackrel{f}{\dashrightarrow}$ & ${\textrm{Y}}$
\end{tabular}

\end{center} It follows that ${\textrm{X}}$ is a weak Del Pezzo surface of degree $9-r$ and $D=3H-Q_1-\ldots -Q_r$.
\begin{example}\label{ex:f1_exm_degree_one}\textrm{\textbf{(degree one Del Pezzo with four cusps)}}
   Let $(L,r)=(\{$ $1123$, $1345$, $1156$, $1258$, $1367$, $1247$, $1468$, $1178$ $~\},5)$ with Dynkin type $4A_2$ in Table~\ref{tab:f1_C1} at row index 163.

 Suppose we want to plant $n$ trees so that there are $3$ trees in each row. The Orchards planting problem is to find a configuration of trees such that the number of rows is maximal. For $n=9$ the solution of this problem is called the Pappus configuration, which is depicted in the image below. We find that there are $10$ rows, where a row is represented by a line.
\begin{center}
 {\includegraphics[width=3cm,height=3cm]{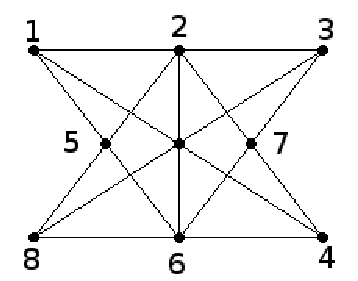}}
\end{center} The labeled points represent 8 base points in the plane. For example the C1 label element $1123$ denotes that a line goes through the base points labeled $1,2$ and $3$. A pair of two lines that do not intersect in the base points form an $A_2$ singularity (for example $1123$ and $1468$). We find three pairs of lines in the Pappus configuration. The fourth pair, $1178$ and $1345$, are no real lines. It is not possible to choose all points real.

 It is well known (see for example \cite{jha2}, chapter 2, section 4, page 281) that a line through two flexes of a smooth cubic curve passes through a third one. So we can choose the base points on eight of the nine flexes of a cubic plane curve.
\end{example}

\section{Algorithm for real structures of weak Del Pezzo surfaces}
A \textit{real surface}\index{real surface} $({\textrm{X}},\sigma)$ is defined as a nonsingular complex projective surface ${\textrm{X}}$ and a complex conjugation $\ensuremath{{\textrm{X}}\stackrel{\sigma}{\rightarrow}{\textrm{X}}}$. See \cite{man1} and \cite{sil1} for theory on real surfaces.

 We define a \textit{real structure}\index{real structure} $\ensuremath{A({\textrm{X}})\stackrel{\sigma_*}{\rightarrow}A({\textrm{X}})}$ of $({\textrm{X}},\sigma)$ as the convolution of its enhanced Picard group $A({\textrm{X}})$ that is uniquely defined by $\sigma$.

 In this section we present an algorithm for classifying the real structures of weak Del Pezzo surfaces up to conjugacy. The classification itself is known and can be found in \cite{wal1}. The correctness of our algorithm is based on the fact that the conjugacy classes of real structures are uniquely determined by the Weyl equivalence classes of root subsystems.
\begin{algorithm}\label{alg:f1_real_alg}\textbf{\textrm{(real structures of weak Del Pezzo surfaces)}}{\small

 \textbf{Output:} Returns a list of unique representatives of equivalence classes of all real structures for weak Del Pezzo surfaces. The list consists of entries $(L,M)$ where $M$ is a matrix representing a real structure, and $(L,$rank$(M))$ is a C1 label that represents the corresponding root subsystem.
\begin{itemize}\addtolength{\itemsep}{1pt}
\item[$\bullet$] In this algorithm a C1 label $(L,r)$ will also be represented as a set of column vectors of length $r+1$  with respect to  the basis ${\textrm{Z}}\langle H,Q_1,\ldots, Q_r\rangle $. For example $1123$ represents the column vector $(1,-1,-1,-1,0,\ldots,0)$. We will represent the equivalence class of a real structure by a matrix acting on these column vectors.
\item[$\bullet$] $C1:=$getGeometricC1RootSubsystemsClassification() (see Algorithm~\ref{alg:f1_geo_basis})
\item[$\bullet$] Let $T$ be an empty ordered list.
\item[$\bullet$] For all $(L,rank)\in C1$  such that  $L \subset \{1123,12,23,34,45,56,67,78\}$ we do the following 4 steps.
\begin{itemize}\addtolength{\itemsep}{1pt}
\item[$\bullet$] $V:=$ a matrix forming a basis consisting of the column vectors $L$ appended by column vectors that are orthogonal to $L$  with respect to  intersection product. See Example~\ref{ex:f1_extend}.
\item[$\bullet$] $D:=$ a square diagonal matrix with rank equal to the rank of $V$. The first $n$ diagonal entries are $-1$, and $1$ otherwise.
\item[$\bullet$] $M:=V\cdot D\cdot V^{-1}$
\item[$\bullet$] If the matrix $M$ has only integer entries then append $(L,M)$ to $T$.
\end{itemize}
\item[$\bullet$] Return $T$.
\end{itemize} }
\end{algorithm}

\begin{example}\label{ex:f1_extend}\textrm{\textbf{(extending to a basis)}}
  Our implementation for constructing $V$ from $L$ in Algorithm~\ref{alg:f1_real_alg} can be best explained by the following 3 examples. The reason is that pseudo code would involve some trickery with indices. Let $V_i$ be the matrix resulting from $L_i$ where $L_1=\{$ $34$, $45$ $\}$, $L_2=\{$ $23$, $34$, $45$ $\}$ and $L_3=\{$ $1123$, $12$, $23$, $45$ $\}$. These correspond to the C1 labels with rank 5 in Table~\ref{tab:f1_C1} at  respectively  the indices 20, 23 and 26. In the matrices below we add $|$ for illustrating the separation between the input columns from $L$ and the extended columns. \[{\tiny V_1= \begin{bmatrix}
  0  &  0  &|&  1  &  0  &  0  &  0\\
  0  &  0  &|&  0  &  1  &  0  &  0\\
  0  &  0  &|&  0  &  0  &  1  &  0\\
  1  &  0  &|&  0  &  0  &  0  &  1\\
 -1  &  1  &|&  0  &  0  &  0  &  1\\
  0  & -1  &|&  0  &  0  &  0  &  1\\

\end{bmatrix}
 V_2= \begin{bmatrix}
  0  &  0  &  0  &|& 1  &  0  &  0\\
  0  &  0  &  0  &|& 0  &  1  &  0\\
  1  &  0  &  0  &|& 0  &  0  &  1\\
 -1  &  1  &  0  &|& 0  &  0  &  1\\
  0  & -1  &  1  &|& 0  &  0  &  1\\
  0  &  0  & -1  &|& 0  &  0  &  1\\

\end{bmatrix}
 V_3= \begin{bmatrix}
  1 &  0 &  0 &  0 &|& -3 & 0 \\
 -1 &  1 &  0 &  0 &|&  1 & 0 \\
 -1 & -1 &  1 &  0 &|&  1 & 0 \\
 -1 &  0 & -1 &  0 &|&  1 & 0 \\
  0 &  0 &  0 &  1 &|&  0 & 1 \\
  0 &  0 &  0 & -1 &|&  0 & 1 \\

\end{bmatrix}
 } \] For each row with zeros for the first 2 columns of $V_1$ we add a canonical basis column vector. For each sequence $Q_i-Q_{i+1}$, $Q_{i+1}-Q_{i+2}$, $\ldots$, $Q_{i+j}-Q_{i+j+1}$ we add a column vector that has a $1$ at indices in the range $[i,j+1]$ and is $0$ otherwise. The same is illustrated for $V_2$. The situation for $V_3$ is a bit different since it contains a column represented by $1123$. For such a column we add the transpose of $(-3,1,1,1,0,0)$. The we continue the same game as in $V_1$ and $V_2$ for the matrix without the first 4 rows.

 By contruction the column vectors are orthogonal. It is a combinatorial exercise to show that the column vectors are linear independent. Finally one has to show that the above described procedure leads to a full rank matrix. The easiest way to do this is by simple checking this fact on all the C1 labels of Table~\ref{tab:f1_C1}.

\end{example}

\begin{example}\label{ex:f1_real_alg}\textrm{\textbf{(real structures of weak Del Pezzo surfaces)}}
  Let $(L,r)=(\{$ $1123$, $12$, $23$, $45$ $\},5)$ be a C1 label of type $A_2+2A_1$ (index 26 in Table~\ref{tab:f1_C1}).  We trace Algorithm~\ref{alg:f1_real_alg} line by line for this C1 label.  We have that $V$ equals $V_3$ in Example~\ref{ex:f1_extend}.  It follows that \[{\tiny V= \begin{bmatrix}
  1 &  0 &  0 &  0 & -3 & 0 \\
 -1 &  1 &  0 &  0 &  1 & 0 \\
 -1 & -1 &  1 &  0 &  1 & 0 \\
 -1 &  0 & -1 &  0 &  1 & 0 \\
  0 &  0 &  0 &  1 &  0 & 1 \\
  0 &  0 &  0 & -1 &  0 & 1 \\

\end{bmatrix}
 , D= \begin{bmatrix}
 -1 &  0 &  0 &  0 &  0 & 0 \\
  0 & -1 &  0 &  0 &  0 & 0 \\
  0 &  0 & -1 &  0 &  0 & 0 \\
  0 &  0 &  0 & -1 &  0 & 0 \\
  0 &  0 &  0 &  0 &  1 & 0 \\
  0 &  0 &  0 &  0 &  0 & 1 \\

\end{bmatrix}
 , M= \begin{bmatrix}
  2 &     1 &     1 &     1 & 0 & 0 \\
 -1 &  -\frac{4}{3} &  -\frac{1}{3} &  -\frac{1}{3} & 0 & 0 \\
 -1 &  -\frac{1}{3} &  -\frac{4}{3} &  -\frac{1}{3} & 0 & 0 \\
 -1 &  -\frac{1}{3} &  -\frac{1}{3} &  -\frac{4}{3} & 0 & 0 \\
  0 &     0 &     0 &     0 & 0 & 1 \\
  0 &     0 &     0 &     0 & 1 & 0 \\

\end{bmatrix}
 }. \]

 We find that $M$ is not an integral matrix and thus not a valid real structure.

\end{example}

\begin{theorem}\label{thm:f1_real_alg}\textbf{\textrm{(real structures of weak Del Pezzo surfaces)}}

 Algorithm~\ref{alg:f1_real_alg} is correct.

\begin{proof}
 Let $({\textrm{X}},\sigma)$ be a real weak Del Pezzo surface.
 Let $W=A({\textrm{X}})=(~W,~K,~\cdot,~h~)$ be its enhanced Picard group.
 Let $D$, $M$ and $V$ be the matrices as in Algorithm~\ref{alg:f1_real_alg}.
 See Example~\ref{ex:f1_extend} for the correctness of the specification of $V$.
 From the real structure $\ensuremath{W\stackrel{\sigma_*}{\rightarrow}W}$ being an involution it follows that $W=W^+\oplus W^-$, where $W^+=\{w\in W|\sigma(w)=w\}$ and $W^-=\{w\in W|\sigma(w)=-w\}$. From $\sigma_*(K)=K$ and $CK=\sigma_*(C)\sigma_*(K)=-CK$ for $C\in W^-$ it follows that $W^-$ is orthogonal to $K$. We have that $W^-$ is generated by a subset of a basis for the root system $R$. It follows that $W^-$ is generated by a basis of a root subsystem of $R$  such that  the basis is a subset of $\{1123,12,23,34,45,56,67,78\}$. Up to conjugacy we only have to consider the root subsystems of the C1 classification. The eigenvectors of the eigenvalue $-1$ are determined by a C1 label. The hyperplanes orthogonal to the roots in $W^-$ are point wise fixed. The eigenvectors of the eigenvalue $1$ are linear independent vectors that are orthogonal to eigenspace of $-1$. It follows that $MV=VD$ and thus $M=VDV^{-1}$. Since $M$ should induce a ${\textbf{Z}}$-module automorphism we need that the matrix $M$ is integral. \end{proof}
\end{theorem}
The output of Algorithm~\ref{alg:f1_real_alg} is in Table~\ref{tab:f1_real}. From Table~\ref{tab:f1_real} it follows that $(rank,f0,f1)$ is a strong invariant for the real structure.

\section{Classification of families of real conics on cyclides}

\label{sec:f1_cyclide}In \textsection\ref{sec:f1_def} we defined the enhanced Picard group and the Weyl object of a weak Del Pezzo surface. We considered weak Del Pezzo surfaces equivalent if their enhanced Picard groups are isomorphic. From Proposition~\ref{prop:f1_weyl} it follows that enhanced Picard groups are isomorphic  if and only if  their Weyl objects are isomorphic as root subsystems. The latter are classified by Algorithm~\ref{alg:f1_geo_basis}.

 The \textit{real enhanced Picard group}\index{real enhanced Picard group} $A({\textrm{X}},\sigma)$ is defined as $ (~\textrm{Pic}({\textrm{X}}),~K,~\cdot,~ h~,~\sigma_*~) $ where $\ensuremath{A({\textrm{X}})\stackrel{\sigma_*}{\rightarrow}A({\textrm{X}})}$ is the real structure. An isomorphism of real enhanced Picard groups is an isomorphism of Picard groups which is compatible with the real structure.

 In this section we introduce a method for classifying real enhanced Picard groups up to isomorphism, by means of a case study.

 We define the \textit{anticanonical model}\index{anticanonical model} of a weak Del Pezzo surface ${\textrm{X}}$ as the image of ${\textrm{X}}$ in projective space under the map associated to the global sections $H^0(-K)$ of the anticanonical divisor $-K$. An \textit{algebraic family of curves}\index{algebraic family of curves} of ${\textrm{X}}$ indexed by a smooth algebraic curve ${\textrm{I}}$ is defined as an irreducible codimension one algebraic subset $F\subset{\textrm{I}}\times{\textrm{X}}$. In this section we assume that the curves are conics  with respect to  the anticanonical model of ${\textrm{X}}$ in projective space.

 We consider weak Del Pezzo surfaces of degree 4. The anticanonical model of such Del Pezzo surfaces lives in projective 4-space. The projection of such anticanonical model to 3-space from an outside point is classically called a \textit{cyclide}\index{cyclide} (see section 8.6.2 in \cite{dol1}). As an application we classify families of real conics on cyclides. First we show that families of conics on cyclides are uniquely defined by certain divisor classes.
\begin{proposition}\label{prop:f1_T2}\textbf{\textrm{(divisor class of a family of conics)}}

 Let ${\textrm{X}}$ a weak Del Pezzo surface of degree $2<K^2<9$.

 A family of conics of the anticanonical model of ${\textrm{X}}$ is uniquely defined by a divisor class in the (2,0)-set that has a positive intersection product with all the effective classes in the (0,-2)-set.

\begin{proof} If $K^2>2$ then from theorem 8.3.2 in \cite{dol1} it follows that the anticanonical model  with respect to  $H^0(-K)$ is a surface. If $K^2<9$ then there exists conics on the anticanonical model.

 Let $C'$ be a class in the (2,0)-set of ${\textrm{X}}$.

 Using Riemann-Roch theorem (see for example \cite{mat1}) it follows $h^0(C')-h^1(C')+h^2(C')=\frac{1}{2}C'(C'-K)+1=2$. From Serre duality, $-K$ being nef and $-K(K-C')<0$ it follows that $h^2(C')=h^0(K-C')=0$. Thus $h^0(C')>1$.

 Let $C'=C+S$ be the decomposition of $C'$ in the mobile part $C$ and the fixed part $S$. Suppose by contradiction that $-KC=-KS=1$. We have that divisors with classes $C$ or $S$ are lines on the anticanonical model. From the adjunction formula it follows that $C^2+KC=-2$ and thus $C^2=1$. Similarly it follows that $S^2=1$. We have that $(C+S)^2=0$ and $CS\geq 0$.Contradiction. From $-K$ being nef and big it follows that $-KC>0$. It follows that $KS=0$ and $-KC=2$. From the proof of lemma 8.2.18 in \cite{dol1} it follows that a fixed class $S$  such that  $KS=0$ is a sum of effective classes in the (0,-2)-set. Moreover, we have that $C=C'$  if and only if  $C'$ is positive against the effective classes in the (0,-2)-set.

 From the vanishing theorem (see chapter 4 in \cite{laz1}) and $-K+C$ being nef and big it follows that $h^0(C)=2$. Now we can define a family of conics $F\subset{\textbf{P}}^1\times{\textrm{X}}$ by the fibers of the morphism $\ensuremath{{\textrm{X}}\stackrel{\varphi_C}{\rightarrow}{\textbf{P}}^1}$ associated to the global sections $H^0(C)$. Thus the curves in the linear series $|C|$ defines $F$ and the family $F$ defines a complete linear series.

 Conversely, suppose that $C$ is the divisor class of a conic in a family of conics of ${\textrm{X}}$. From the adjunction formula it follows that $C^2+CK=-2$. From $-KC=2$ it follows that $C^2=0$ and thus $C$ is in the (2,0)-set. By definition $C$ has no fixed components and is thus positive against the effective classes in the (0,-2)-set. \end{proof}
\end{proposition}
Let $({\textrm{X}},\sigma)$ be a real weak Del Pezzo surface of degree 4.
\begin{itemize}\addtolength{\itemsep}{1pt}
\item[$\bullet$] Let $C$ be the set of C1 labels of root subsystems of $D_5$ (indices $16,\ldots,31$ in Table~\ref{tab:f1_C1}).
\item[$\bullet$] Let $R$ be the set of conjugacy classes of real structures of ${\textrm{X}}$ (indices $10,\ldots,15$ in Table~\ref{tab:f1_real}).
\item[$\bullet$] Let $B$ be the set of subsets of the roots in the (0,-2)-set of ${\textrm{X}}$ that have a geometric C1 label.
\item[$\bullet$] Let $B(c)=\{ b\in B | b {\textrm{~and~}} c \textrm{ are isomorphic bases} \}$ for $c\in C$ (see Proposition~\ref{prop:f1_C_K}).
\item[$\bullet$] Let $B(c,r)=\{ b\in B(c) | r(b)=b \}$ for $r\in R$ (note that we need that the basis of a root subsystem is compatible with the real structure).
\item[$\bullet$] Let $G(c,r,b)$ be the set of classes in the (2,0)-set that are fixed by the real structure $r\in R$ and have positive intersection product against classes $b$ in the real effective (0,-2)-set $B(c,r)$.
\end{itemize} The following table denotes the possible cardinalities of $G(c,r,b)$ by fixing a representative $r$ for a conjugacy class of a real structure in $R$, and considering all bases in $b\in B(c)$ of the same Dynkin type as $c$.

 The `index' column denotes the index of $c$ in Table~\ref{tab:f1_C1} and the `type' column, the Dynkin type of $c$. The `index' row denote the index $i$ for the real structure $r$ in Table~\ref{tab:f1_real} and the `type' row, the Dynkin type of $r$.
\begin{center}
{\tiny \begin{tabular}
{|c|c||c|c|c|c|c|c|} \hline index &          & $10$ & $11$  & $12$ & $13$  & $14$ & $15$ \\
\hline       & type     & $A0$ & $A1$  & $2A1$& $2A1$ & $3A1$& $D4$ \\
\hline \hline $16$ & $ A0    $ & $10$ & $6$   & $2$  & $6  $ & $2 $ & $2 $ \\
\hline $17$ & $ A1    $ & $8 $ & $4,6$ & $2$  & $4  $ & $2 $ & $  $ \\
\hline $18$ & $2A1    $ & $6 $ & $4$   & $2$  & $   $ & $  $ & $  $ \\
\hline $19$ & $2A1    $ & $7 $ & $3$   & $ $  & $3,5$ & $1 $ & $1 $ \\
\hline $20$ & $ A2    $ & $6 $ & $2$   & $ $  & $2  $ & $2 $ & $  $ \\
\hline $21$ & $3A1    $ & $5 $ & $3$   & $ $  & $3  $ & $1 $ & $  $ \\
\hline $22$ & $ A2+ A1$ & $4 $ & $2$   & $ $  & $   $ & $  $ & $  $ \\
\hline $23$ & $ A3    $ & $4 $ & $ $   & $2$  & $   $ & $  $ & $  $ \\
\hline $24$ & $ A3    $ & $5 $ & $1$   & $ $  & $1,3$ & $  $ & $  $ \\
\hline $25$ & $4A1    $ & $4 $ & $ $   & $2$  & $2,4$ & $  $ & $  $ \\
\hline $26$ & $ A2+2A1$ & $3 $ & $ $   & $ $  & $   $ & $1 $ & $  $ \\
\hline $27$ & $ A3+ A1$ & $3 $ & $1$   & $ $  & $   $ & $  $ & $  $ \\
\hline $28$ & $ A4    $ & $2 $ & $ $   & $ $  & $   $ & $  $ & $  $ \\
\hline $29$ & $ D4    $ & $3 $ & $ $   & $ $  & $1  $ & $  $ & $  $ \\
\hline $30$ & $ A3+2A1$ & $2 $ & $ $   & $ $  & $2  $ & $  $ & $  $ \\
\hline $31$ & $ D5    $ & $1 $ & $ $   & $ $  & $   $ & $  $ & $  $ \\
\hline
\end{tabular}
 }
\end{center} For example the anticanonical model of a real torus in 4-space has 4 complex $A_1$ singularities and 4 families of real conics. The singularities and the real structure of the torus are defined by entry $(row,column)=(25,13)$ in the above table. Two of the families of conics are defined by the Villarceau circles.

 We see from the table that the number of families of conics on ${\textrm{X}}$ is not uniquely defined by only the conjugacy classes of $c$ and $r$. For example if $(row,column)=(17,11)$ then there exists $b=[12],b'=[1123]\in B(c)$  such that  $\#G(c,r,b)=4$ and $\#G(c,r,b')=6$. We have that $r: (\alpha_0,\alpha_1,\alpha_2,\alpha_3,\alpha_4,\alpha_5)=\alpha_0H+\ensuremath{\overset{}{\underset{}{\sum}}}\alpha_iQ_i \mapsto (\alpha_0,\alpha_1,\alpha_2,\alpha_3,\alpha_5,\alpha_4)$, $G(c,r,b)=\{$ $H-Q_1$, $H-Q_3$, $2H-Q_1-Q_2-Q_4-Q_5$, $2H-Q_1-Q_3-Q_4-Q_5$ $\}$ and $G(c,r,b')=\{$ $H-Q_1$, $H-Q_2$, $H-Q_3$, $2H-Q_1-Q_2-Q_4-Q_5$, $2H-Q_1-Q_3-Q_4-Q_5$, $2H-Q_2-Q_3-Q_4-Q_5$ $\}$.

 The obvious algorithm computing the above table does not terminate in reasonable time for degree 1 and 2 weak Del Pezzo surfaces, since $B$ is too large. The roots in the eigenspace of $-1$ cannot be effective. It follows that we only have to consider bases $b\in B$ that are disjoint to the eigenspace of $-1$.

\section{Acknowledgements}
It is my pleasure to acknowledge that the many computations with Josef Schicho is a major contribution to this paper.
Also he recognized the Pappus configuration of example 7.
I would like to thank Michael Harrison for useful discussions concerning root subsystems.
I would like to thank Ulrich Derenthal for informing me of a mistake in a previous verion of this paper.

The algorithms were implemented using the computer algebra system Sage (\cite{sage}).

This research was supported by the Austrian Science Fund (FWF): project P21461.\appendix

\section{Table of geometric C1 labels}

\begin{tab}\label{tab:f1_C1}\textbf{\textrm{(geometric bases of C1 root subsystems up to Weyl equivalence)}}
\begin{itemize}\addtolength{\itemsep}{1pt}
\item[$\bullet$] the `index' column assigns an index to each row for future reference,
\item[$\bullet$] the `rank' column defines the rank $r$ of the C1 root system,
\item[$\bullet$] the `basis' column denotes the basis $L$  such that  $(L,r)$ is a geometric C1 label,
\item[$\bullet$] the `type' column denotes the Dynkin type of the C1 label $(L,r)$.
\end{itemize}
{\tiny
\begin{longtable}
{| c|c|c|c| } \hline index & rank & basis & type  \\
 \hline \hline $ 1   $ & $ 2   $ & $ []                                               $ & $  A0             $ \\
 \hline $ 2   $ & $ 2   $ & $ [12]                                             $ & $  A1             $ \\
 \hline \hline $ 3   $ & $ 3   $ & $ []                                               $ & $  A0             $ \\
 \hline $ 4   $ & $ 3   $ & $ [23]                                             $ & $  A1             $ \\
 \hline $ 5   $ & $ 3   $ & $ [1123]                                           $ & $  A1             $ \\
 \hline $ 6   $ & $ 3   $ & $ [1123, 23]                                       $ & $ 2A1             $ \\
 \hline $ 7   $ & $ 3   $ & $ [12, 23]                                         $ & $  A2             $ \\
 \hline $ 8   $ & $ 3   $ & $ [1123, 12, 23]                                   $ & $  A2+ A1         $ \\
 \hline \hline $ 9   $ & $ 4   $ & $ []                                               $ & $  A0             $ \\
 \hline $ 10  $ & $ 4   $ & $ [34]                                             $ & $  A1             $ \\
 \hline $ 11  $ & $ 4   $ & $ [12, 34]                                         $ & $ 2A1             $ \\
 \hline $ 12  $ & $ 4   $ & $ [23, 34]                                         $ & $  A2             $ \\
 \hline $ 13  $ & $ 4   $ & $ [1123, 12, 34]                                   $ & $  A2+ A1         $ \\
 \hline $ 14  $ & $ 4   $ & $ [12, 23, 34]                                     $ & $  A3             $ \\
 \hline $ 15  $ & $ 4   $ & $ [1123, 12, 23, 34]                               $ & $  A4             $ \\
 \hline \hline $ 16  $ & $ 5   $ & $ []                                               $ & $  A0             $ \\
 \hline $ 17  $ & $ 5   $ & $ [45]                                             $ & $  A1             $ \\
 \hline $ 18  $ & $ 5   $ & $ [23, 45]                                         $ & $ 2A1             $ \\
 \hline $ 19  $ & $ 5   $ & $ [1123, 45]                                       $ & $ 2A1             $ \\
 \hline $ 20  $ & $ 5   $ & $ [34, 45]                                         $ & $  A2             $ \\
 \hline $ 21  $ & $ 5   $ & $ [1123, 23, 45]                                   $ & $ 3A1             $ \\
 \hline $ 22  $ & $ 5   $ & $ [12, 34, 45]                                     $ & $  A2+ A1         $ \\
 \hline $ 23  $ & $ 5   $ & $ [23, 34, 45]                                     $ & $  A3             $ \\
 \hline $ 24  $ & $ 5   $ & $ [1123, 34, 45]                                   $ & $  A3             $ \\
 \hline $ 25  $ & $ 5   $ & $ [1145, 1123, 23, 45]                             $ & $ 4A1             $ \\
 \hline $ 26  $ & $ 5   $ & $ [1123, 12, 23, 45]                               $ & $  A2+2A1         $ \\
 \hline $ 27  $ & $ 5   $ & $ [1123, 12, 34, 45]                               $ & $  A3+ A1         $ \\
 \hline $ 28  $ & $ 5   $ & $ [12, 23, 34, 45]                                 $ & $  A4             $ \\
 \hline $ 29  $ & $ 5   $ & $ [1123, 23, 34, 45]                               $ & $  D4             $ \\
 \hline $ 30  $ & $ 5   $ & $ [1145, 1123, 12, 23, 45]                         $ & $  A3+2A1         $ \\
 \hline $ 31  $ & $ 5   $ & $ [1123, 12, 23, 34, 45]                           $ & $  D5             $ \\
 \hline \hline $ 32  $ & $ 6   $ & $ []                                               $ & $  A0             $ \\
 \hline $ 33  $ & $ 6   $ & $ [56]                                             $ & $  A1             $ \\
 \hline $ 34  $ & $ 6   $ & $ [34, 56]                                         $ & $ 2A1             $ \\
 \hline $ 35  $ & $ 6   $ & $ [45, 56]                                         $ & $  A2             $ \\
 \hline $ 36  $ & $ 6   $ & $ [12, 34, 56]                                     $ & $ 3A1             $ \\
 \hline $ 37  $ & $ 6   $ & $ [23, 45, 56]                                     $ & $  A2+ A1         $ \\
 \hline $ 38  $ & $ 6   $ & $ [34, 45, 56]                                     $ & $  A3             $ \\
 \hline $ 39  $ & $ 6   $ & $ [1145, 1123, 23, 45]                             $ & $ 4A1             $ \\
 \hline $ 40  $ & $ 6   $ & $ [1123, 23, 45, 56]                               $ & $  A2+2A1         $ \\
 \hline $ 41  $ & $ 6   $ & $ [12, 23, 45, 56]                                 $ & $ 2A2             $ \\
 \hline $ 42  $ & $ 6   $ & $ [12, 34, 45, 56]                                 $ & $  A3+ A1         $ \\
 \hline $ 43  $ & $ 6   $ & $ [23, 34, 45, 56]                                 $ & $  A4             $ \\
 \hline $ 44  $ & $ 6   $ & $ [1123, 23, 34, 45]                               $ & $  D4             $ \\
 \hline $ 45  $ & $ 6   $ & $ [1123, 12, 23, 45, 56]                           $ & $ 2A2+ A1         $ \\
 \hline $ 46  $ & $ 6   $ & $ [1145, 1123, 23, 45, 56]                         $ & $  A3+2A1         $ \\
 \hline $ 47  $ & $ 6   $ & $ [1123, 12, 34, 45, 56]                           $ & $  A4+ A1         $ \\
 \hline $ 48  $ & $ 6   $ & $ [12, 23, 34, 45, 56]                             $ & $  A5             $ \\
 \hline $ 49  $ & $ 6   $ & $ [1123, 23, 34, 45, 56]                           $ & $  D5             $ \\
 \hline $ 50  $ & $ 6   $ & $ [1456, 1123, 12, 23, 45, 56]                     $ & $ 3A2             $ \\
 \hline $ 51  $ & $ 6   $ & $ [1145, 1123, 12, 23, 45, 56]                     $ & $  A5+ A1         $ \\
 \hline $ 52  $ & $ 6   $ & $ [1123, 12, 23, 34, 45, 56]                       $ & $  E6             $ \\
 \hline \hline $ 53  $ & $ 7   $ & $ []                                               $ & $  A0             $ \\
 \hline $ 54  $ & $ 7   $ & $ [67]                                             $ & $  A1             $ \\
 \hline $ 55  $ & $ 7   $ & $ [45, 67]                                         $ & $ 2A1             $ \\
 \hline $ 56  $ & $ 7   $ & $ [56, 67]                                         $ & $  A2             $ \\
 \hline $ 57  $ & $ 7   $ & $ [23, 45, 67]                                     $ & $ 3A1             $ \\
 \hline $ 58  $ & $ 7   $ & $ [1123, 45, 67]                                   $ & $ 3A1             $ \\
 \hline $ 59  $ & $ 7   $ & $ [34, 56, 67]                                     $ & $  A2+ A1         $ \\
 \hline $ 60  $ & $ 7   $ & $ [45, 56, 67]                                     $ & $  A3             $ \\
 \hline $ 61  $ & $ 7   $ & $ [1123, 23, 45, 67]                               $ & $ 4A1             $ \\
 \hline $ 62  $ & $ 7   $ & $ [1145, 1123, 23, 45]                             $ & $ 4A1             $ \\
 \hline $ 63  $ & $ 7   $ & $ [12, 34, 56, 67]                                 $ & $  A2+2A1         $ \\
 \hline $ 64  $ & $ 7   $ & $ [23, 34, 56, 67]                                 $ & $ 2A2             $ \\
 \hline $ 65  $ & $ 7   $ & $ [23, 45, 56, 67]                                 $ & $  A3+ A1         $ \\
 \hline $ 66  $ & $ 7   $ & $ [1123, 45, 56, 67]                               $ & $  A3+ A1         $ \\
 \hline $ 67  $ & $ 7   $ & $ [34, 45, 56, 67]                                 $ & $  A4             $ \\
 \hline $ 68  $ & $ 7   $ & $ [1123, 23, 34, 45]                               $ & $  D4             $ \\
 \hline $ 69  $ & $ 7   $ & $ [1145, 1123, 23, 45, 67]                         $ & $ 5A1             $ \\
 \hline $ 70  $ & $ 7   $ & $ [1123, 12, 23, 45, 67]                           $ & $  A2+3A1         $ \\
 \hline $ 71  $ & $ 7   $ & $ [1123, 12, 34, 56, 67]                           $ & $ 2A2+ A1         $ \\
 \hline $ 72  $ & $ 7   $ & $ [1123, 23, 45, 56, 67]                           $ & $  A3+2A1         $ \\
 \hline $ 73  $ & $ 7   $ & $ [1145, 1123, 23, 45, 56]                         $ & $  A3+2A1         $ \\
 \hline $ 74  $ & $ 7   $ & $ [12, 23, 45, 56, 67]                             $ & $  A3+ A2         $ \\
 \hline $ 75  $ & $ 7   $ & $ [12, 34, 45, 56, 67]                             $ & $  A4+ A1         $ \\
 \hline $ 76  $ & $ 7   $ & $ [23, 34, 45, 56, 67]                             $ & $  A5             $ \\
 \hline $ 77  $ & $ 7   $ & $ [1123, 34, 45, 56, 67]                           $ & $  A5             $ \\
 \hline $ 78  $ & $ 7   $ & $ [1123, 23, 34, 45, 67]                           $ & $  D4+ A1         $ \\
 \hline $ 79  $ & $ 7   $ & $ [1123, 23, 34, 45, 56]                           $ & $  D5             $ \\
 \hline $ 80  $ & $ 7   $ & $ [1567, 1347, 1127, 12, 34, 56]                   $ & $ 6A1             $ \\
 \hline $ 81  $ & $ 7   $ & $ [1456, 1123, 12, 23, 45, 56]                     $ & $ 3A2             $ \\
 \hline $ 82  $ & $ 7   $ & $ [1145, 1123, 12, 23, 45, 67]                     $ & $  A3+3A1         $ \\
 \hline $ 83  $ & $ 7   $ & $ [1123, 12, 23, 45, 56, 67]                       $ & $  A3+ A2+ A1     $ \\
 \hline $ 84  $ & $ 7   $ & $ [278, 12, 23, 34, 56, 67]                        $ & $ 2A3             $ \\
 \hline $ 85  $ & $ 7   $ & $ [1123, 12, 23, 34, 56, 67]                       $ & $  A4+ A2         $ \\
 \hline $ 86  $ & $ 7   $ & $ [1123, 12, 34, 45, 56, 67]                       $ & $  A5+ A1         $ \\
 \hline $ 87  $ & $ 7   $ & $ [1145, 1123, 12, 23, 45, 56]                     $ & $  A5+ A1         $ \\
 \hline $ 88  $ & $ 7   $ & $ [12, 23, 34, 45, 56, 67]                         $ & $  A6             $ \\
 \hline $ 89  $ & $ 7   $ & $ [1145, 1123, 23, 45, 56, 67]                     $ & $  D4+2A1         $ \\
 \hline $ 90  $ & $ 7   $ & $ [1123, 12, 23, 34, 45, 67]                       $ & $  D5+ A1         $ \\
 \hline $ 91  $ & $ 7   $ & $ [1123, 23, 34, 45, 56, 67]                       $ & $  D6             $ \\
 \hline $ 92  $ & $ 7   $ & $ [1123, 12, 23, 34, 45, 56]                       $ & $  E6             $ \\
 \hline $ 93  $ & $ 7   $ & $ [278, 1567, 1347, 1127, 12, 34, 56]              $ & $ 7A1             $ \\
 \hline $ 94  $ & $ 7   $ & $ [278, 1567, 12, 23, 34, 56, 67]                  $ & $ 2A3+ A1         $ \\
 \hline $ 95  $ & $ 7   $ & $ [1456, 1123, 12, 23, 45, 56, 67]                 $ & $  A5+ A2         $ \\
 \hline $ 96  $ & $ 7   $ & $ [278, 12, 23, 34, 45, 56, 67]                    $ & $  A7             $ \\
 \hline $ 97  $ & $ 7   $ & $ [278, 1347, 1127, 12, 34, 45, 56]                $ & $  D4+3A1         $ \\
 \hline $ 98  $ & $ 7   $ & $ [1145, 1123, 12, 23, 45, 56, 67]                 $ & $  D6+ A1         $ \\
 \hline $ 99  $ & $ 7   $ & $ [1123, 12, 23, 34, 45, 56, 67]                   $ & $  E7             $ \\
 \hline \hline $ 100 $ & $ 8   $ & $ []                                               $ & $  A0             $ \\
 \hline $ 101 $ & $ 8   $ & $ [78]                                             $ & $  A1             $ \\
 \hline $ 102 $ & $ 8   $ & $ [56, 78]                                         $ & $ 2A1             $ \\
 \hline $ 103 $ & $ 8   $ & $ [67, 78]                                         $ & $  A2             $ \\
 \hline $ 104 $ & $ 8   $ & $ [34, 56, 78]                                     $ & $ 3A1             $ \\
 \hline $ 105 $ & $ 8   $ & $ [45, 67, 78]                                     $ & $  A2+ A1         $ \\
 \hline $ 106 $ & $ 8   $ & $ [56, 67, 78]                                     $ & $  A3             $ \\
 \hline $ 107 $ & $ 8   $ & $ [12, 34, 56, 78]                                 $ & $ 4A1             $ \\
 \hline $ 108 $ & $ 8   $ & $ [1145, 1123, 23, 45]                             $ & $ 4A1             $ \\
 \hline $ 109 $ & $ 8   $ & $ [23, 45, 67, 78]                                 $ & $  A2+2A1         $ \\
 \hline $ 110 $ & $ 8   $ & $ [34, 45, 67, 78]                                 $ & $ 2A2             $ \\
 \hline $ 111 $ & $ 8   $ & $ [34, 56, 67, 78]                                 $ & $  A3+ A1         $ \\
 \hline $ 112 $ & $ 8   $ & $ [45, 56, 67, 78]                                 $ & $  A4             $ \\
 \hline $ 113 $ & $ 8   $ & $ [1123, 23, 34, 45]                               $ & $  D4             $ \\
 \hline $ 114 $ & $ 8   $ & $ [1145, 1123, 23, 45, 78]                         $ & $ 5A1             $ \\
 \hline $ 115 $ & $ 8   $ & $ [1123, 23, 45, 67, 78]                           $ & $  A2+3A1         $ \\
 \hline $ 116 $ & $ 8   $ & $ [12, 34, 45, 67, 78]                             $ & $ 2A2+ A1         $ \\
 \hline $ 117 $ & $ 8   $ & $ [12, 34, 56, 67, 78]                             $ & $  A3+2A1         $ \\
 \hline $ 118 $ & $ 8   $ & $ [1145, 1123, 23, 45, 56]                         $ & $  A3+2A1         $ \\
 \hline $ 119 $ & $ 8   $ & $ [23, 34, 56, 67, 78]                             $ & $  A3+ A2         $ \\
 \hline $ 120 $ & $ 8   $ & $ [23, 45, 56, 67, 78]                             $ & $  A4+ A1         $ \\
 \hline $ 121 $ & $ 8   $ & $ [34, 45, 56, 67, 78]                             $ & $  A5             $ \\
 \hline $ 122 $ & $ 8   $ & $ [1123, 23, 34, 45, 78]                           $ & $  D4+ A1         $ \\
 \hline $ 123 $ & $ 8   $ & $ [1123, 23, 34, 45, 56]                           $ & $  D5             $ \\
 \hline $ 124 $ & $ 8   $ & $ [1567, 1347, 1127, 12, 34, 56]                   $ & $ 6A1             $ \\
 \hline $ 125 $ & $ 8   $ & $ [1145, 1123, 23, 45, 67, 78]                     $ & $  A2+4A1         $ \\
 \hline $ 126 $ & $ 8   $ & $ [1123, 12, 23, 45, 67, 78]                       $ & $ 2A2+2A1         $ \\
 \hline $ 127 $ & $ 8   $ & $ [1456, 1123, 12, 23, 45, 56]                     $ & $ 3A2             $ \\
 \hline $ 128 $ & $ 8   $ & $ [1145, 1123, 23, 45, 56, 78]                     $ & $  A3+3A1         $ \\
 \hline $ 129 $ & $ 8   $ & $ [1123, 12, 34, 56, 67, 78]                       $ & $  A3+ A2+ A1     $ \\
 \hline $ 130 $ & $ 8   $ & $ [12, 23, 34, 56, 67, 78]                         $ & $ 2A3             $ \\
 \hline $ 131 $ & $ 8   $ & $ [1567, 1145, 1127, 12, 56, 78]                   $ & $ 2A3             $ \\
 \hline $ 132 $ & $ 8   $ & $ [1123, 23, 45, 56, 67, 78]                       $ & $  A4+2A1         $ \\
 \hline $ 133 $ & $ 8   $ & $ [12, 23, 45, 56, 67, 78]                         $ & $  A4+ A2         $ \\
 \hline $ 134 $ & $ 8   $ & $ [12, 34, 45, 56, 67, 78]                         $ & $  A5+ A1         $ \\
 \hline $ 135 $ & $ 8   $ & $ [1145, 1123, 12, 23, 45, 56]                     $ & $  A5+ A1         $ \\
 \hline $ 136 $ & $ 8   $ & $ [23, 34, 45, 56, 67, 78]                         $ & $  A6             $ \\
 \hline $ 137 $ & $ 8   $ & $ [1145, 1123, 23, 45, 56, 67]                     $ & $  D4+2A1         $ \\
 \hline $ 138 $ & $ 8   $ & $ [1123, 23, 34, 45, 67, 78]                       $ & $  D4+ A2         $ \\
 \hline $ 139 $ & $ 8   $ & $ [1123, 23, 34, 45, 56, 78]                       $ & $  D5+ A1         $ \\
 \hline $ 140 $ & $ 8   $ & $ [1123, 23, 34, 45, 56, 67]                       $ & $  D6             $ \\
 \hline $ 141 $ & $ 8   $ & $ [1123, 12, 23, 34, 45, 56]                       $ & $  E6             $ \\
 \hline $ 142 $ & $ 8   $ & $ [278, 1567, 1347, 1127, 12, 34, 56]              $ & $ 7A1             $ \\
 \hline $ 143 $ & $ 8   $ & $ [1456, 1123, 12, 23, 45, 56, 78]                 $ & $ 3A2+ A1         $ \\
 \hline $ 144 $ & $ 8   $ & $ [278, 1347, 1127, 12, 34, 56, 78]                $ & $  A3+4A1         $ \\
 \hline $ 145 $ & $ 8   $ & $ [1145, 1123, 12, 23, 45, 67, 78]                 $ & $  A3+ A2+2A1     $ \\
 \hline $ 146 $ & $ 8   $ & $ [1678, 1456, 1347, 1123, 12, 45, 78]             $ & $ 2A3+ A1         $ \\
 \hline $ 147 $ & $ 8   $ & $ [1123, 12, 23, 45, 56, 67, 78]                   $ & $  A4+ A2+ A1     $ \\
 \hline $ 148 $ & $ 8   $ & $ [1123, 12, 23, 34, 56, 67, 78]                   $ & $  A4+ A3         $ \\
 \hline $ 149 $ & $ 8   $ & $ [1145, 1123, 12, 23, 45, 56, 78]                 $ & $  A5+2A1         $ \\
 \hline $ 150 $ & $ 8   $ & $ [1456, 1123, 12, 23, 45, 56, 67]                 $ & $  A5+ A2         $ \\
 \hline $ 151 $ & $ 8   $ & $ [1123, 12, 34, 45, 56, 67, 78]                   $ & $  A6+ A1         $ \\
 \hline $ 152 $ & $ 8   $ & $ [12, 23, 34, 45, 56, 67, 78]                     $ & $  A7             $ \\
 \hline $ 153 $ & $ 8   $ & $ [1567, 1145, 1127, 12, 23, 56, 78]               $ & $  A7             $ \\
 \hline $ 154 $ & $ 8   $ & $ [1567, 1347, 1127, 12, 34, 56, 78]               $ & $  D4+3A1         $ \\
 \hline $ 155 $ & $ 8   $ & $ [1567, 1347, 1145, 1127, 12, 56, 78]             $ & $  D4+ A3         $ \\
 \hline $ 156 $ & $ 8   $ & $ [1145, 1123, 23, 45, 56, 67, 78]                 $ & $  D5+2A1         $ \\
 \hline $ 157 $ & $ 8   $ & $ [1123, 12, 23, 34, 45, 67, 78]                   $ & $  D5+ A2         $ \\
 \hline $ 158 $ & $ 8   $ & $ [1145, 1123, 12, 23, 45, 56, 67]                 $ & $  D6+ A1         $ \\
 \hline $ 159 $ & $ 8   $ & $ [1123, 23, 34, 45, 56, 67, 78]                   $ & $  D7             $ \\
 \hline $ 160 $ & $ 8   $ & $ [1123, 12, 23, 34, 45, 56, 78]                   $ & $  E6+ A1         $ \\
 \hline $ 161 $ & $ 8   $ & $ [1123, 12, 23, 34, 45, 56, 67]                   $ & $  E7             $ \\
 \hline $ 162 $ & $ 8   $ & $ [308, 278, 1567, 1347, 1127, 12, 34, 56]         $ & $ 8A1             $ \\
 \hline $ 163 $ & $ 8   $ & $ [1123, 1345, 1156, 1258, 1367, 1247, 1468, 1178] $ & $ 4A2             $ \\
 \hline $ 164 $ & $ 8   $ & $ [308, 278, 1567, 12, 23, 34, 56, 67]             $ & $ 2A3+2A1         $ \\
 \hline $ 165 $ & $ 8   $ & $ [278, 1678, 12, 23, 34, 45, 67, 78]              $ & $ 2A4             $ \\
 \hline $ 166 $ & $ 8   $ & $ [1678, 1145, 1123, 12, 23, 45, 67, 78]           $ & $  A5+ A2+ A1     $ \\
 \hline $ 167 $ & $ 8   $ & $ [234, 1145, 1123, 12, 23, 56, 67, 78]            $ & $  A7+ A1         $ \\
 \hline $ 168 $ & $ 8   $ & $ [1567, 1123, 12, 23, 34, 56, 67, 78]             $ & $  A8             $ \\
 \hline $ 169 $ & $ 8   $ & $ [278, 1567, 1347, 1127, 12, 34, 56, 78]          $ & $  D4+4A1         $ \\
 \hline $ 170 $ & $ 8   $ & $ [234, 278, 12, 23, 34, 56, 67, 78]               $ & $ 2D4             $ \\
 \hline $ 171 $ & $ 8   $ & $ [278, 1567, 12, 23, 34, 56, 67, 78]              $ & $  D5+ A3         $ \\
 \hline $ 172 $ & $ 8   $ & $ [278, 1347, 1127, 12, 34, 45, 56, 78]            $ & $  D6+2A1         $ \\
 \hline $ 173 $ & $ 8   $ & $ [278, 12, 23, 34, 45, 56, 67, 78]                $ & $  D8             $ \\
 \hline $ 174 $ & $ 8   $ & $ [1456, 1123, 12, 23, 45, 56, 67, 78]             $ & $  E6+ A2         $ \\
 \hline $ 175 $ & $ 8   $ & $ [1145, 1123, 12, 23, 45, 56, 67, 78]             $ & $  E7+ A1         $ \\
 \hline $ 176 $ & $ 8   $ & $ [1123, 12, 23, 34, 45, 56, 67, 78]               $ & $  E8             $ \\
 \hline
\end{longtable}
}
\end{tab}

\section{Table of real structures of weak Del Pezzo surfaces}
The following table represents the output of Algorithm~\ref{alg:f1_real_alg}.
\begin{tab}\label{tab:f1_real}\textbf{\textrm{(real structures of weak Del Pezzo surfaces)}}
\begin{itemize}\addtolength{\itemsep}{1pt}
\item[$\bullet$] the `index' column denotes the row index for future reference,
\item[$\bullet$] the `rank'  column denotes that the root subsystem of the real structure is contained in a root system of given rank,
\item[$\bullet$] the `C1'    column denotes the index in Table~\ref{tab:f1_C1} corresponding to the root subsystem of the real structure,
\item[$\bullet$] the `type'  column denotes the type of the root root subsystem of the real structure,
\item[$\bullet$] the `f0' columns denote the cardinality of the real (0,-2)-set (thus classes in the (0,-2)-set that are fixed by the real structure),
\item[$\bullet$] the `f1' columns denote the cardinality of the real (1,-1)-set, and
\item[$\bullet$] the `f2' columns denote the cardinality of the real (2,0)-set.
\end{itemize}
{\tiny
\begin{longtable}
{|c|c|c|c| c|c|c|c|} \hline index & rank & C1 & type & f0 & f1 & f2 \\
 \hline \hline $ 1 $ & $ 2 $ & $ 1 $ & $  A0 $ & $ 2 $ & $ 3 $ & $ 2 $ \\
 \hline $ 2 $ & $ 2 $ & $ 2 $ & $  A1 $ & $ 0 $ & $ 1 $ & $ 0 $ \\
 \hline \hline $ 3 $ & $ 3 $ & $ 3 $ & $  A0 $ & $ 8 $ & $ 6 $ & $ 3 $ \\
 \hline $ 4 $ & $ 3 $ & $ 4 $ & $  A1 $ & $ 2 $ & $ 2 $ & $ 1 $ \\
 \hline $ 5 $ & $ 3 $ & $ 5 $ & $  A1 $ & $ 6 $ & $ 0 $ & $ 3 $ \\
 \hline $ 6 $ & $ 3 $ & $ 6 $ & $ 2A1 $ & $ 0 $ & $ 0 $ & $ 1 $ \\
 \hline \hline $ 7 $ & $ 4 $ & $ 9 $ & $  A0 $ & $ 20 $ & $ 10 $ & $ 5 $ \\
 \hline $ 8 $ & $ 4 $ & $ 10 $ & $  A1 $ & $ 6 $ & $ 4 $ & $ 3 $  \\
 \hline $ 9 $ & $ 4 $ & $ 11 $ & $ 2A1 $ & $ 0 $ & $ 2 $ & $ 1 $  \\
 \hline \hline $ 10 $ & $ 5 $ & $ 16 $ & $  A0 $ & $ 40 $ & $ 16 $ & $ 10 $ \\
 \hline $ 11 $ & $ 5 $ & $ 17 $ & $  A1 $ & $ 14 $ & $ 8 $ & $ 6 $   \\
 \hline $ 12 $ & $ 5 $ & $ 18 $ & $ 2A1 $ & $ 4 $ & $ 4 $ & $ 2 $    \\
 \hline $ 13 $ & $ 5 $ & $ 19 $ & $ 2A1 $ & $ 12 $ & $ 0 $ & $ 6 $   \\
 \hline $ 14 $ & $ 5 $ & $ 21 $ & $ 3A1 $ & $ 2 $ & $ 0 $ & $ 2 $    \\
 \hline $ 15 $ & $ 5 $ & $ 29 $ & $  D4 $ & $ 0 $ & $ 0 $ & $ 2 $    \\
 \hline \hline $ 16 $ & $ 6 $ & $ 32 $ & $  A0 $ & $ 72 $ & $ 27 $ & $ 27 $        \\
 \hline $ 17 $ & $ 6 $ & $ 33 $ & $  A1 $ & $ 30 $ & $ 15 $ & $ 15 $        \\
 \hline $ 18 $ & $ 6 $ & $ 34 $ & $ 2A1 $ & $ 12 $ & $ 7 $ & $ 7 $          \\
 \hline $ 19 $ & $ 6 $ & $ 36 $ & $ 3A1 $ & $ 2 $ & $ 3 $ & $ 3 $           \\
 \hline $ 20 $ & $ 6 $ & $ 44 $ & $  D4 $ & $ 0 $ & $ 3 $ & $ 3 $           \\
 \hline \hline $ 21 $ & $ 7 $ & $ 53 $ & $  A0 $ & $ 126 $ & $ 56 $ & $ 126 $      \\
 \hline $ 22 $ & $ 7 $ & $ 54 $ & $  A1 $ & $ 60 $ & $ 32 $ & $ 60 $        \\
 \hline $ 23 $ & $ 7 $ & $ 55 $ & $ 2A1 $ & $ 26 $ & $ 16 $ & $ 26 $        \\
 \hline $ 24 $ & $ 7 $ & $ 57 $ & $ 3A1 $ & $ 8 $ & $ 8 $ & $ 8 $           \\
 \hline $ 25 $ & $ 7 $ & $ 58 $ & $ 3A1 $ & $ 24 $ & $ 0 $ & $ 24 $         \\
 \hline $ 26 $ & $ 7 $ & $ 61 $ & $ 4A1 $ & $ 6 $ & $ 0 $ & $ 6 $           \\
 \hline $ 27 $ & $ 7 $ & $ 68 $ & $  D4 $ & $ 6 $ & $ 8 $ & $ 6 $           \\
 \hline $ 28 $ & $ 7 $ & $ 78 $ & $  D4+ A1 $ & $ 4 $ & $ 0 $ & $ 4 $       \\
 \hline $ 29 $ & $ 7 $ & $ 91 $ & $  D6 $ & $ 2 $ & $ 0 $ & $ 2 $           \\
 \hline $ 30 $ & $ 7 $ & $ 99 $ & $  E7 $ & $ 0 $ & $ 0 $ & $ 0 $           \\
 \hline \hline $ 31 $ & $ 8 $ & $ 100 $ & $  A0 $ & $ 240 $ & $ 240 $ & $ 2160 $   \\
 \hline $ 32 $ & $ 8 $ & $ 101 $ & $  A1 $ & $ 126 $ & $ 126 $ & $ 756 $    \\
 \hline $ 33 $ & $ 8 $ & $ 102 $ & $ 2A1 $ & $ 60 $ & $ 60 $ & $ 252 $      \\
 \hline $ 34 $ & $ 8 $ & $ 104 $ & $ 3A1 $ & $ 26 $ & $ 26 $ & $ 72 $       \\
 \hline $ 35 $ & $ 8 $ & $ 107 $ & $ 4A1 $ & $ 8 $ & $ 8 $ & $ 24 $         \\
 \hline $ 36 $ & $ 8 $ & $ 113 $ & $  D4 $ & $ 24 $ & $ 24 $ & $ 24 $       \\
 \hline $ 37 $ & $ 8 $ & $ 122 $ & $  D4+ A1 $ & $ 6 $ & $ 6 $ & $ 12 $     \\
 \hline $ 38 $ & $ 8 $ & $ 140 $ & $  D6 $ & $ 4 $ & $ 4 $ & $ 4 $          \\
 \hline $ 39 $ & $ 8 $ & $ 161 $ & $  E7 $ & $ 2 $ & $ 2 $ & $ 0 $          \\
 \hline $ 40 $ & $ 8 $ & $ 176 $ & $  E8 $ & $ 0 $ & $ 0 $ & $ 0 $          \\
 \hline
\end{longtable}
 }
Below is an explicit coordinate description of the real structure for each index in the above table: \[ \ensuremath{\sigma_*:A({\textrm{X}})\rightarrow A({\textrm{X}}),\quad (H,Q_1,\ldots,Q_r)\mapsto (D_0,\ldots,D_r)} \] where $A({\textrm{X}})$ is the enhanced Picard group of a real weak Del Pezzo surface $({\textrm{X}},\sigma)$. Below we denote $(D_0,\ldots,D_r)$ for each index.{\tiny  $\bullet$ $index = 1 $,  ( $H$  $,$  $Q_1$  $,$  $Q_2$  )  $\bullet$ $index = 2 $,  ( $H$  $,$  $Q_2$  $,$  $Q_1$  )  $\bullet$ $index = 3 $,  ( $H$  $,$  $Q_1$  $,$  $Q_2$  $,$  $Q_3$  )  $\bullet$ $index = 4 $,  ( $H$  $,$  $Q_1$  $,$  $Q_3$  $,$  $Q_2$  )  $\bullet$ $index = 5 $,  ( $2H$ $-$ $Q_1$ $-$ $Q_2$ $-$ $Q_3$  $,$  $H$ $-$ $Q_2$ $-$ $Q_3$  $,$  $H$ $-$ $Q_1$ $-$ $Q_3$  $,$  $H$ $-$ $Q_1$ $-$ $Q_2$  )  $\bullet$ $index = 6 $,  ( $2H$ $-$ $Q_1$ $-$ $Q_2$ $-$ $Q_3$  $,$  $H$ $-$ $Q_2$ $-$ $Q_3$  $,$  $H$ $-$ $Q_1$ $-$ $Q_2$  $,$  $H$ $-$ $Q_1$ $-$ $Q_3$  )  $\bullet$ $index = 7 $,  ( $H$  $,$  $Q_1$  $,$  $Q_2$  $,$  $Q_3$  $,$  $Q_4$  )  $\bullet$ $index = 8 $,  ( $H$  $,$  $Q_1$  $,$  $Q_2$  $,$  $Q_4$  $,$  $Q_3$  )  $\bullet$ $index = 9 $,  ( $H$  $,$  $Q_2$  $,$  $Q_1$  $,$  $Q_4$  $,$  $Q_3$  )  $\bullet$ $index = 10 $,  ( $H$  $,$  $Q_1$  $,$  $Q_2$  $,$  $Q_3$  $,$  $Q_4$  $,$  $Q_5$  )  $\bullet$ $index = 11 $,  ( $H$  $,$  $Q_1$  $,$  $Q_2$  $,$  $Q_3$  $,$  $Q_5$  $,$  $Q_4$  )  $\bullet$ $index = 12 $,  ( $H$  $,$  $Q_1$  $,$  $Q_3$  $,$  $Q_2$  $,$  $Q_5$  $,$  $Q_4$  )  $\bullet$ $index = 13 $,  ( $2H$ $-$ $Q_1$ $-$ $Q_2$ $-$ $Q_3$  $,$  $H$ $-$ $Q_2$ $-$ $Q_3$  $,$  $H$ $-$ $Q_1$ $-$ $Q_3$  $,$  $H$ $-$ $Q_1$ $-$ $Q_2$  $,$  $Q_5$  $,$  $Q_4$  )  $\bullet$ $index = 14 $,  ( $2H$ $-$ $Q_1$ $-$ $Q_2$ $-$ $Q_3$  $,$  $H$ $-$ $Q_2$ $-$ $Q_3$  $,$  $H$ $-$ $Q_1$ $-$ $Q_2$  $,$  $H$ $-$ $Q_1$ $-$ $Q_3$  $,$  $Q_5$  $,$  $Q_4$  )  $\bullet$ $index = 15 $,  ( $3H$ $-$ $2Q_1$ $-$ $Q_2$ $-$ $Q_3$ $-$ $Q_4$ $-$ $Q_5$  $,$  $2H$ $-$ $Q_1$ $-$ $Q_2$ $-$ $Q_3$ $-$ $Q_4$ $-$ $Q_5$  $,$  $H$ $-$ $Q_1$ $-$ $Q_2$  $,$  $H$ $-$ $Q_1$ $-$ $Q_3$  $,$  $H$ $-$ $Q_1$ $-$ $Q_4$  $,$  $H$ $-$ $Q_1$ $-$ $Q_5$  )  $\bullet$ $index = 16 $,  ( $H$  $,$  $Q_1$  $,$  $Q_2$  $,$  $Q_3$  $,$  $Q_4$  $,$  $Q_5$  $,$  $Q_6$  )  $\bullet$ $index = 17 $,  ( $H$  $,$  $Q_1$  $,$  $Q_2$  $,$  $Q_3$  $,$  $Q_4$  $,$  $Q_6$  $,$  $Q_5$  )  $\bullet$ $index = 18 $,  ( $H$  $,$  $Q_1$  $,$  $Q_2$  $,$  $Q_4$  $,$  $Q_3$  $,$  $Q_6$  $,$  $Q_5$  )  $\bullet$ $index = 19 $,  ( $H$  $,$  $Q_2$  $,$  $Q_1$  $,$  $Q_4$  $,$  $Q_3$  $,$  $Q_6$  $,$  $Q_5$  )  $\bullet$ $index = 20 $,  ( $3H$ $-$ $2Q_1$ $-$ $Q_2$ $-$ $Q_3$ $-$ $Q_4$ $-$ $Q_5$  $,$  $2H$ $-$ $Q_1$ $-$ $Q_2$ $-$ $Q_3$ $-$ $Q_4$ $-$ $Q_5$  $,$  $H$ $-$ $Q_1$ $-$ $Q_2$  $,$  $H$ $-$ $Q_1$ $-$ $Q_3$  $,$  $H$ $-$ $Q_1$ $-$ $Q_4$  $,$  $H$ $-$ $Q_1$ $-$ $Q_5$  $,$  $Q_6$  )  $\bullet$ $index = 21 $,  ( $H$  $,$  $Q_1$  $,$  $Q_2$  $,$  $Q_3$  $,$  $Q_4$  $,$  $Q_5$  $,$  $Q_6$  $,$  $Q_7$  )  $\bullet$ $index = 22 $,  ( $H$  $,$  $Q_1$  $,$  $Q_2$  $,$  $Q_3$  $,$  $Q_4$  $,$  $Q_5$  $,$  $Q_7$  $,$  $Q_6$  )  $\bullet$ $index = 23 $,  ( $H$  $,$  $Q_1$  $,$  $Q_2$  $,$  $Q_3$  $,$  $Q_5$  $,$  $Q_4$  $,$  $Q_7$  $,$  $Q_6$  )  $\bullet$ $index = 24 $,  ( $H$  $,$  $Q_1$  $,$  $Q_3$  $,$  $Q_2$  $,$  $Q_5$  $,$  $Q_4$  $,$  $Q_7$  $,$  $Q_6$  )  $\bullet$ $index = 25 $,  ( $2H$ $-$ $Q_1$ $-$ $Q_2$ $-$ $Q_3$  $,$  $H$ $-$ $Q_2$ $-$ $Q_3$  $,$  $H$ $-$ $Q_1$ $-$ $Q_3$  $,$  $H$ $-$ $Q_1$ $-$ $Q_2$  $,$  $Q_5$  $,$  $Q_4$  $,$  $Q_7$  $,$  $Q_6$  )  $\bullet$ $index = 26 $,  ( $2H$ $-$ $Q_1$ $-$ $Q_2$ $-$ $Q_3$  $,$  $H$ $-$ $Q_2$ $-$ $Q_3$  $,$  $H$ $-$ $Q_1$ $-$ $Q_2$  $,$  $H$ $-$ $Q_1$ $-$ $Q_3$  $,$  $Q_5$  $,$  $Q_4$  $,$  $Q_7$  $,$  $Q_6$  )  $\bullet$ $index = 27 $,  ( $3H$ $-$ $2Q_1$ $-$ $Q_2$ $-$ $Q_3$ $-$ $Q_4$ $-$ $Q_5$  $,$  $2H$ $-$ $Q_1$ $-$ $Q_2$ $-$ $Q_3$ $-$ $Q_4$ $-$ $Q_5$  $,$  $H$ $-$ $Q_1$ $-$ $Q_2$  $,$  $H$ $-$ $Q_1$ $-$ $Q_3$  $,$  $H$ $-$ $Q_1$ $-$ $Q_4$  $,$  $H$ $-$ $Q_1$ $-$ $Q_5$  $,$  $Q_6$  $,$  $Q_7$  )  $\bullet$ $index = 28 $,  ( $3H$ $-$ $2Q_1$ $-$ $Q_2$ $-$ $Q_3$ $-$ $Q_4$ $-$ $Q_5$  $,$  $2H$ $-$ $Q_1$ $-$ $Q_2$ $-$ $Q_3$ $-$ $Q_4$ $-$ $Q_5$  $,$  $H$ $-$ $Q_1$ $-$ $Q_2$  $,$  $H$ $-$ $Q_1$ $-$ $Q_3$  $,$  $H$ $-$ $Q_1$ $-$ $Q_4$  $,$  $H$ $-$ $Q_1$ $-$ $Q_5$  $,$  $Q_7$  $,$  $Q_6$  )  $\bullet$ $index = 29 $,  ( $4H$ $-$ $3Q_1$ $-$ $Q_2$ $-$ $Q_3$ $-$ $Q_4$ $-$ $Q_5$ $-$ $Q_6$ $-$ $Q_7$  $,$  $3H$ $-$ $2Q_1$ $-$ $Q_2$ $-$ $Q_3$ $-$ $Q_4$ $-$ $Q_5$ $-$ $Q_6$ $-$ $Q_7$  $,$  $H$ $-$ $Q_1$ $-$ $Q_2$  $,$  $H$ $-$ $Q_1$ $-$ $Q_3$  $,$  $H$ $-$ $Q_1$ $-$ $Q_4$  $,$  $H$ $-$ $Q_1$ $-$ $Q_5$  $,$  $H$ $-$ $Q_1$ $-$ $Q_6$  $,$  $H$ $-$ $Q_1$ $-$ $Q_7$  )  $\bullet$ $index = 30 $,  ( $8H$ $-$ $3Q_1$ $-$ $3Q_2$ $-$ $3Q_3$ $-$ $3Q_4$ $-$ $3Q_5$ $-$ $3Q_6$ $-$ $3Q_7$  $,$  $3H$ $-$ $2Q_1$ $-$ $Q_2$ $-$ $Q_3$ $-$ $Q_4$ $-$ $Q_5$ $-$ $Q_6$ $-$ $Q_7$  $,$  $3H$ $-$ $Q_1$ $-$ $2Q_2$ $-$ $Q_3$ $-$ $Q_4$ $-$ $Q_5$ $-$ $Q_6$ $-$ $Q_7$  $,$  $3H$ $-$ $Q_1$ $-$ $Q_2$ $-$ $2Q_3$ $-$ $Q_4$ $-$ $Q_5$ $-$ $Q_6$ $-$ $Q_7$  $,$  $3H$ $-$ $Q_1$ $-$ $Q_2$ $-$ $Q_3$ $-$ $2Q_4$ $-$ $Q_5$ $-$ $Q_6$ $-$ $Q_7$  $,$  $3H$ $-$ $Q_1$ $-$ $Q_2$ $-$ $Q_3$ $-$ $Q_4$ $-$ $2Q_5$ $-$ $Q_6$ $-$ $Q_7$  $,$  $3H$ $-$ $Q_1$ $-$ $Q_2$ $-$ $Q_3$ $-$ $Q_4$ $-$ $Q_5$ $-$ $2Q_6$ $-$ $Q_7$  $,$  $3H$ $-$ $Q_1$ $-$ $Q_2$ $-$ $Q_3$ $-$ $Q_4$ $-$ $Q_5$ $-$ $Q_6$ $-$ $2Q_7$  )  $\bullet$ $index = 31 $,  ( $H$  $,$  $Q_1$  $,$  $Q_2$  $,$  $Q_3$  $,$  $Q_4$  $,$  $Q_5$  $,$  $Q_6$  $,$  $Q_7$  $,$  $Q_8$  )  $\bullet$ $index = 32 $,  ( $H$  $,$  $Q_1$  $,$  $Q_2$  $,$  $Q_3$  $,$  $Q_4$  $,$  $Q_5$  $,$  $Q_6$  $,$  $Q_8$  $,$  $Q_7$  )  $\bullet$ $index = 33 $,  ( $H$  $,$  $Q_1$  $,$  $Q_2$  $,$  $Q_3$  $,$  $Q_4$  $,$  $Q_6$  $,$  $Q_5$  $,$  $Q_8$  $,$  $Q_7$  )  $\bullet$ $index = 34 $,  ( $H$  $,$  $Q_1$  $,$  $Q_2$  $,$  $Q_4$  $,$  $Q_3$  $,$  $Q_6$  $,$  $Q_5$  $,$  $Q_8$  $,$  $Q_7$  )  $\bullet$ $index = 35 $,  ( $H$  $,$  $Q_2$  $,$  $Q_1$  $,$  $Q_4$  $,$  $Q_3$  $,$  $Q_6$  $,$  $Q_5$  $,$  $Q_8$  $,$  $Q_7$  )  $\bullet$ $index = 36 $,  ( $3H$ $-$ $2Q_1$ $-$ $Q_2$ $-$ $Q_3$ $-$ $Q_4$ $-$ $Q_5$  $,$  $2H$ $-$ $Q_1$ $-$ $Q_2$ $-$ $Q_3$ $-$ $Q_4$ $-$ $Q_5$  $,$  $H$ $-$ $Q_1$ $-$ $Q_2$  $,$  $H$ $-$ $Q_1$ $-$ $Q_3$  $,$  $H$ $-$ $Q_1$ $-$ $Q_4$  $,$  $H$ $-$ $Q_1$ $-$ $Q_5$  $,$  $Q_6$  $,$  $Q_7$  $,$  $Q_8$  )  $\bullet$ $index = 37 $,  ( $3H$ $-$ $2Q_1$ $-$ $Q_2$ $-$ $Q_3$ $-$ $Q_4$ $-$ $Q_5$  $,$  $2H$ $-$ $Q_1$ $-$ $Q_2$ $-$ $Q_3$ $-$ $Q_4$ $-$ $Q_5$  $,$  $H$ $-$ $Q_1$ $-$ $Q_2$  $,$  $H$ $-$ $Q_1$ $-$ $Q_3$  $,$  $H$ $-$ $Q_1$ $-$ $Q_4$  $,$  $H$ $-$ $Q_1$ $-$ $Q_5$  $,$  $Q_6$  $,$  $Q_8$  $,$  $Q_7$  )  $\bullet$ $index = 38 $,  ( $4H$ $-$ $3Q_1$ $-$ $Q_2$ $-$ $Q_3$ $-$ $Q_4$ $-$ $Q_5$ $-$ $Q_6$ $-$ $Q_7$  $,$  $3H$ $-$ $2Q_1$ $-$ $Q_2$ $-$ $Q_3$ $-$ $Q_4$ $-$ $Q_5$ $-$ $Q_6$ $-$ $Q_7$  $,$  $H$ $-$ $Q_1$ $-$ $Q_2$  $,$  $H$ $-$ $Q_1$ $-$ $Q_3$  $,$  $H$ $-$ $Q_1$ $-$ $Q_4$  $,$  $H$ $-$ $Q_1$ $-$ $Q_5$  $,$  $H$ $-$ $Q_1$ $-$ $Q_6$  $,$  $H$ $-$ $Q_1$ $-$ $Q_7$  $,$  $Q_8$  )  $\bullet$ $index = 39 $,  ( $8H$ $-$ $3Q_1$ $-$ $3Q_2$ $-$ $3Q_3$ $-$ $3Q_4$ $-$ $3Q_5$ $-$ $3Q_6$ $-$ $3Q_7$  $,$  $3H$ $-$ $2Q_1$ $-$ $Q_2$ $-$ $Q_3$ $-$ $Q_4$ $-$ $Q_5$ $-$ $Q_6$ $-$ $Q_7$  $,$  $3H$ $-$ $Q_1$ $-$ $2Q_2$ $-$ $Q_3$ $-$ $Q_4$ $-$ $Q_5$ $-$ $Q_6$ $-$ $Q_7$  $,$  $3H$ $-$ $Q_1$ $-$ $Q_2$ $-$ $2Q_3$ $-$ $Q_4$ $-$ $Q_5$ $-$ $Q_6$ $-$ $Q_7$  $,$  $3H$ $-$ $Q_1$ $-$ $Q_2$ $-$ $Q_3$ $-$ $2Q_4$ $-$ $Q_5$ $-$ $Q_6$ $-$ $Q_7$  $,$  $3H$ $-$ $Q_1$ $-$ $Q_2$ $-$ $Q_3$ $-$ $Q_4$ $-$ $2Q_5$ $-$ $Q_6$ $-$ $Q_7$  $,$  $3H$ $-$ $Q_1$ $-$ $Q_2$ $-$ $Q_3$ $-$ $Q_4$ $-$ $Q_5$ $-$ $2Q_6$ $-$ $Q_7$  $,$  $3H$ $-$ $Q_1$ $-$ $Q_2$ $-$ $Q_3$ $-$ $Q_4$ $-$ $Q_5$ $-$ $Q_6$ $-$ $2Q_7$  $,$  $Q_8$  )  $\bullet$ $index = 40 $,  ( $17H$ $-$ $6Q_1$ $-$ $6Q_2$ $-$ $6Q_3$ $-$ $6Q_4$ $-$ $6Q_5$ $-$ $6Q_6$ $-$ $6Q_7$ $-$ $6Q_8$  $,$  $6H$ $-$ $3Q_1$ $-$ $2Q_2$ $-$ $2Q_3$ $-$ $2Q_4$ $-$ $2Q_5$ $-$ $2Q_6$ $-$ $2Q_7$ $-$ $2Q_8$  $,$  $6H$ $-$ $2Q_1$ $-$ $3Q_2$ $-$ $2Q_3$ $-$ $2Q_4$ $-$ $2Q_5$ $-$ $2Q_6$ $-$ $2Q_7$ $-$ $2Q_8$  $,$  $6H$ $-$ $2Q_1$ $-$ $2Q_2$ $-$ $3Q_3$ $-$ $2Q_4$ $-$ $2Q_5$ $-$ $2Q_6$ $-$ $2Q_7$ $-$ $2Q_8$  $,$  $6H$ $-$ $2Q_1$ $-$ $2Q_2$ $-$ $2Q_3$ $-$ $3Q_4$ $-$ $2Q_5$ $-$ $2Q_6$ $-$ $2Q_7$ $-$ $2Q_8$  $,$  $6H$ $-$ $2Q_1$ $-$ $2Q_2$ $-$ $2Q_3$ $-$ $2Q_4$ $-$ $3Q_5$ $-$ $2Q_6$ $-$ $2Q_7$ $-$ $2Q_8$  $,$  $6H$ $-$ $2Q_1$ $-$ $2Q_2$ $-$ $2Q_3$ $-$ $2Q_4$ $-$ $2Q_5$ $-$ $3Q_6$ $-$ $2Q_7$ $-$ $2Q_8$  $,$  $6H$ $-$ $2Q_1$ $-$ $2Q_2$ $-$ $2Q_3$ $-$ $2Q_4$ $-$ $2Q_5$ $-$ $2Q_6$ $-$ $3Q_7$ $-$ $2Q_8$  $,$  $6H$ $-$ $2Q_1$ $-$ $2Q_2$ $-$ $2Q_3$ $-$ $2Q_4$ $-$ $2Q_5$ $-$ $2Q_6$ $-$ $2Q_7$ $-$ $3Q_8$  ) }
\end{tab}

\bibliography{geometry,schicho,rational_curves}
\paragraph{Address of author:}
~\\
~\\
King Abdullah University of Science and Technology, Thuwal, Kingdom of Saudi Arabia
\\
\textbf{email:} niels.lubbes@gmail.com
\end{document}